\documentclass[11pt]{amsart}
\usepackage{amsmath}
\usepackage{amssymb}
\usepackage{amscd}
\usepackage{latexsym}
\usepackage{graphics, epsfig}
\usepackage{curves}

\theoremstyle{plain}
\newtheorem{thm}{Theorem}[section]
\newtheorem{lem}[thm]{Lemma}

\newtheorem{exmp}{Example}

\newtheorem*{coro}{Corollary}

\theoremstyle{definition}

\def \CPb {\overline{\mathbf{CP}}^{2}}
\def \CP {{\mathbf{CP}}^{2}} 
\def \R {\mathbf{R}}
\def \Z {\mathbf{Z}}

\def \a {\alpha}
\def \b {\beta}

\def \s {\sigma}

\def \bd {\partial}
\def \x {\times}
\def \- {\setminus}

\def \ssw {\text{SW}}
\def \sw {\mathcal{SW}}

\def \DD {\Delta}

\def\fs{\mathfrak{s}}

\title{Construction of exotic smooth structures}

\begin{document}

\author{Anar Akhmedov}
\address{Department of Mathematics\\
University of California at Irvine\\
Irvine, CA 92697-3875}

\email{aakhmedo@math.uci.edu}

\begin{abstract} In this article, we construct infinitely many simply connected, nonsymplectic and pairwise nondiffeomorphic 4-manifolds starting from the elliptic surfaces $E(n)$ and applying the sequence of knot surgery, ordinary blowups and rational blowdown. We also compute the Seiberg-Witten invariants of these manifolds.  
\end{abstract}

\maketitle

\setcounter{section}{-1}

\section{Introduction}

There has been a recent flurry of activity in the discovery of exotic smooth structures on simply-connected $4$-manifolds with small Euler characteristic. In the 2004, Jongil Park \cite{P2} gave the first example of exotic smooth structure on $\CP\#7\,\CPb$ i.e. 4-manifold homeomorphic to $\CP\#\, 7\CP$ but not diffeomophic to it. Shortly afterwards, Andr\'as Stipsicz and Zolt\'an Szab\'o used a technique similar to Park's to construct an exotic smooth structure on  $\CP\#\, 6\CPb$ \cite{SS1}. Then Fintushel and Stern \cite{FS3} introduced a new technique, surgery in double nodes, which demonstrated that in fact $\CP\#k\,\CPb$, $k= 6,7,8$, have infinitely many distinct smooth structures. Park, Stipsicz, and Szab\'o  \cite{PSS}, using \cite{FS3}, constructed infinitely many smooth structures when $k=5$. Stipsicz and Szab\'o used similar  ideas to construct exotic smooth structures on $3\,\CP\#k\,\CPb$ for $k = 9$ \cite{SS2} and Park for $k = 8$ \cite{P3}. All these infinite family of manifolds were constructed from the elliptic surfaces $E(1)$ and $E(2)$ by applying the sequence of knot surgery in double nodes, oridinary blowups and rational blowdown. In this article, we obtain similar result starting from $E(n)$ for $n\geq 3$. Applying the sequence of knot surgery, ordinary blowups and rational blowdown we construct an infinite family of simply connected, nonsymplectic and pairwise nondiffeomorphic manifolds with nontrivial Seiberg-Witten invariants. Our motivation for constructing such manifolds comes from the following question raised by Ronald Stern.

\medskip

\textit{Question}. Which homotopy types of a simply connected, smooth, closed 4-manifold can be obtained from $E(n)$ by the sequence of knot surgery, ordinary blowups and rational blowdowns?

\medskip

The author was informed by Jongil Park that a similar result has been obtained also by him and Ki-Heon Yun.

{\bf Acknowledgments:} I am grateful to Ron Stern for many helpful discussions. I also wish to thank the referee for the helpful remarks.

\section{Elliptic fibrations}

\subsection{Elliptic fibrations as words in mapping class group }

Let $f\colon X^4\to S^2$ be a genus $1$ Lefschetz fibration. Then $X^4$ is $T^2$ bundle over $S^2$ away from finitely many critical values of $f$. The fibers over these critical values are called singular fibers of the Lefschetz fibration. Such fibrations are characterized by their monodromy, i.e. a factorization of the identity element in the mapping class group $M_{1}$ of the fiber $T^2$ as a product of right-handed Dehn twists. 

It is well-known that the mapping class group $M_{1}$ of the torus is generated by $t_{a}$, $t_{b}$ of $M_{1}$ subject to the relations \[ t_{a}t_{b}t_{a}=t_{b}t_{a}t_{b}\qquad {\mbox{ and }} \qquad (t_{a}t_{b})^6=1 \]

\noindent where $t_{a}$, $t_{b}$ are Dehn twists along the standard curves $a$ and $b$ that generate the first homology of the torus. 

Classification of genus $1$ Lefschetz fibrations were given by Moishezon \cite{M}, who showed that after a possible perturbation such a fibration over $S^2$ is equivalent to one of the fibrations given by the words $(t_{a}t_{b})^{6n}=1$ in $M_{1}$. The total space $X$ of this fibration is $E(n)$, where $E(n)$ is a simply connected elliptic surface. It is known that $E(n)$ has holomorphic Euler characteristic $\chi_{h} = n$ with a section of $E(n)\to S^2$ that is a sphere with self--intersection $-n$.

\subsection{Singular fibers in elliptic fibrations}

In this section we give a brief introduction to certain singular fibers that occur in elliptic fibrations. We also prove the lemma that will be used in our costruction. The complete treatment of the topology of elliptic surfaces and their singular fiberes can be found in \cite{HKK} and \cite{KM}.

\medskip

$\bf {Type \ I_{1}}$ : The fiber is an immersed 2-sphere with one positive double point. It has monodromy that is a conjugate of $t_{a}$ and referred to as a fishtail fiber. 

\medskip

$\bf{Type \ I_{k}}$ : The singularity which has a monodromy conjugate of ${t_{a}}^{k}$ ($k\geq 2$) is called an $I_k$ singularity(or necklace fiber). Singular fibers of type $I_k$ are  a plumbing of $k$ smooth 2--spheres of self--intersection $-2$ along a circle. $I_{2}$ singularity is also called  a double node. 

\medskip

\begin{lem} 

There exists an elliptic Lefschetz fibration on the surface $E(n)$ with a section, a singular fiber $F$ of type $I_{8n}$, $(2n-1)$ singular fibers $F_1, F_2,F_3, ...,F_{2n-1}$ of type $I_2$ and two additional fishtail fibers.

\end{lem}

\begin{proof} The word $(t_{a}t_{b})^{6n}=1$ in the mapping class group $M_1$ determines a genus $1$ Lefschetz fibration on $E(n)$. 
Notice that by using the braid relation $t_{a}t_{b}t_{a} = t_{b}t_{a}t_{b}$, this word becomes equivalent to the following words 

\[ 1 = (t_{a}t_{b})^{6n} = t_{a}(t_{b}t_{a}t_{b})(t_{a}t_{b}t_{a})....(t_{b}t_{a}t_{b})t_{a}t_{b} = \]
\[t_{a}(t_{a}t_{b}t_{a})(t_{b}t_{a}t_{b})....(t_{a}t_{b}t_{a})t_{a}t_{b}= {t_{a}}^{2}(t_{b}t_{a})^{6n-3}t_{b}{t_{a}}^{2}t_{b} =\]
\[{t_{a}}^{2}(t_{b}t_{a}t_{b})(t_{a}t_{b}t_{a})....(t_{a}t_{b}t_{a})t_{b}{t_{a}}^{2}t_{b} = {t_{a}}^{2}(t_{a}t_{b}t_{a})(t_{b}t_{a}t_{b})....(t_{b}t_{a}t_{b})t_{b}{t_{a}}^{2}t_{b} = \]
\[ {t_{a}}^{3}(t_{b}t_{a})^{6n-4}{t_{b}}^{2}{t_{a}}^{2}t_{b} = {t_{a}}^{3}(t_{b}t_{a}t_{b})(t_{a}t_{b}t_{a})....(t_{b}t_{a}t_{b})t_{a}{t_{b}}^{2}{t_{a}}^{2}t_{b} = \]
\[{t_{a}}^{3}(t_{a}t_{b}t_{a})(t_{b}t_{a}t_{b})....(t_{a}t_{b}t_{a})t_{a}{t_{b}}^{2}{t_{a}}^{2}t_{b} = {t_{a}}^{4}(t_{b}t_{a})^{6n-6}t_{b}{t_{a}}^{2}{t_{b}}^{2}{t_{a}}^{2}t_{b} \]

\noindent Applying the braid relation repeteadly many times as above, the word $(t_{a}t_{b})^{6n}=1$  becomes equivalent to 

\[{t_{a}}^{4n}t_{b}{t_{a}}^{2}{t_{b}}^{2}{t_{a}}^{2}{t_{b}}^{2}{t_{a}}^{2}{t_{b}}^{2}.....{t_{a}}^{2}t_{b}=1 \] 

Next, conjugating by powers of $t_{a}$ in the last relation, we have the following word

\[ {t_{a}}^{4n}t_{b}{t_{a}}^{2}{t_{b}}^{2}{t_{a}}^{-2}{t_{a}}^{4}{t_{b}}^{2}{t_{a}}^{-4}.....{t_{a}}^{4n-2}{t_{b}}^{2}{t_{a}}^{-(4n-2)}{t_{a}}^{4n}t_{b}{t_{a}}^{-4n}{t_{a}}^{4n} = 1 \]
 
Finally, multyplying by ${t_{a}}^{4n}$ from the left and ${t_{a}}^{-4n}$ from the right in the last relation, we obtain the following word 

\[ {t_{a}}^{8n}t_{b}{t_{a}}^{2}{t_{b}}^{2}{t_{a}}^{-2}{t_{a}}^{4}{t_{b}}^{2}{t_{a}}^{-4}.....{t_{a}}^{4n-2}{t_{b}}^{2}{t_{a}}^{-(4n-2)}{t_{a}}^{4n}t_{b}{t_{a}}^{-4n}=1 \]
 
\noindent This word translates into the singular fiber of type $I_{8n}$, $(2n-1)$ $I_{2}$ fibers and two fisthtail fibers.   \end{proof}

\section{Seiberg-Witten Invariants}

 In this section we review the basics of Seiberg-Witten invariants introduced by Seiberg and Witten. In order to state theorems in the following sections, we will view the Seiberg-Witten invariant of a smooth $4$-manifold as a multivariable Laurent polynomial \cite {W}. Let us recall that the Seiberg-Witten invariant  of a smooth closed oriented $4$-manifold $X$ with $b_2 ^+(X)>1$ is an integer valued function which is defined on the set of $spin ^{\, c}$ structures over $X$ \cite{W}. For simplicity we assume that $H_1(X,\Z)$ has no 2-torsion. Then there is a one-to-one correspondence between the set of $spin ^{\, c}$ structures over $X$ and the set characteristic elements of $H^2(X,\Z)$ as following: To each $spin ^c$ structure $\fs $ over $X$ corresponds a bundle of positive spinors $W^+_{\fs}$ over $X$. Let $c(\fs)\in H_2(X)$ denote the Poincar\'e dual of $c_1(W^+_{\fs})$. Each $c(\fs)$ is a characteristic element of $H_2(X;\Z)$ (i.e. its Poincar\'e dual $\hat{c}(\fs)=c_1(W^+_{\fs})$ reduces mod~2 to $w_2(X)$).

In this set up we can view the Seiberg-Witten invariant as integer valued function
\[ \ssw_X: \lbrace k\in H^2(X,\Z)|k\equiv w_2(TX)\pmod2)\rbrace
\rightarrow \Z. \] The Seiberg-Witten invariant $\ssw_X$ is a
diffeomorphism invariant and its sign depends on an orientation of
\[ H^0(X,\R)\otimes\det H_+^2(X,\R)\otimes \det H^1(X,\R).\] If
$\ssw_X(\b)\neq 0$, then we call
$\b$ a {\it{basic class}} of $X$. It is a fundamental fact that the set
of basic classes is finite. It can be shown that, if $\b$ is a basic class, then
so is $-\b$ with
\[\ssw_X(-\b)=(-1)^{(\text{e}+\text{sign})(X)/4}\,\ssw_X(\b)\] where
$\text{e}(X)$ is the Euler number and $\text{sign}(X)$ is the signature
of $X$.

Let $\{\pm \b_1,\dots,\pm \b_n\}$ be the set of nonzero basic classes for
$X$. Consider variables $t_{\b}=\exp(\b)$ for each $\b\in H^2(X;\Z)$
which satisfy the relations $t_{\a+\b}=t_{\a}t_{\b}$. We may then view the  Seiberg-Witten invariant of $X$ as the symmetric Laurent polynomial \[\sw_X = b_0+\sum_{j=1}^n b_j(t_{\b_j}+(-1)^{(\text{e}+\text{sign})(X)/4}\, t_{\b_j}^{-1})\] where $b_0=\ssw_X(0)$ and $b_j=\ssw_X(\b_j)$.  

\begin{exmp} Let $E(n)$ be a simply connected minimally elliptic surface with holomorphic Euler characteristic $\chi = n$ and with no multiple fibers. Then we have $\sw_{E(n)}=(t-t^{-1})^{n-2}$ where $t= \exp(T)$ and $T$ is the cohomology class Poincare dual to the fiber class. Thus $\ssw_{E(n)}((n-2i)T)=(-1)^{i-1}\binom{n-2}{i-1}$ for $i=1,\dots,n-1$ and $\ssw_{E(n)}(\b)=0$ for any other $\b$.
\end{exmp}

\medskip

\begin{thm}(Taubes) \cite {T} Suppose that $X$ is a closed symplectic 4-manifold with ${b_{2}}^{+}(X) > 1$. If $K_{X}$ is a canonical classes of $X$, then $SW_{X}(\pm K_{X}) = \pm 1$. \end{thm}

\section{Knot surgery} 

Let $X$ be a $4$-manifold (with ${b_{2}}^{+}(X) > 1 $) which contains a homologically-nontrivial torus $T$ of self-intersection $0$. Let $N(K)$ be a tubular neighborhood of $K$ in $S^3$, and let $T\x D^2$ be a tubular neighborhood of $T$ in $X$. Then the knot surgery manifold $X_K$ is defined by \[ X_K = (X\- (T\x D^2))\cup (S^1\x (S^3\- N(K))\]

\noindent Fintushel and Stern proved the theorem that shows Seiberg-Witten invariants of $X_K$  can be completely determined by the Seiberg-Witten invariant of $X$ and the Alexander polynomial of $K$ \cite{FS2}. Furthermore, if $X$ and $X\- T$ are simply connected, then so is $X_K$.

\begin{thm} \cite {FS2}
Assume that $T$ lies in a cusp neighborhood in $X$, then Seiberg-Witten invariants of $X_{K}$ is $\sw_{X_{K}} = \sw_{X} \centerdot \DD_{K}(t)$. If the Alexander polynomial $\DD_{K}(t)$ of knot $K$ is not monic then $X_{K}$ admits no symplectic structure.

\end{thm}

\noindent We refer the reader to \cite{FS3} for recent technique of Fintushel-Stern on knot surgery in double nodes.

\section{Rational blow-downs and its generalizations} In early nineties, the rational blow-down surgery was introduced by Fintushel and Stern \cite{FS1}. The basic idea is that if a smooth 4-manifold $X$ contains a particular configuration $C_{p}$ of transversally intersecting 2-spheres whose boundary is the lens space $L(p^2,1-p)$, then one can replace $C_{p}$ with rational ball $B_{p}$ to construct a new manifold $X_{p}$. If one knows the Donaldson and the Seiberg-Witten invariants of the original manifold $X$, then Fintushel-Stern show how one can determine same invariants of $X_{p}$. Initially Fintushel and Stern used this surgery technique to compute gauge invariants for simply connected regular elliptic surfaces with multiple fibers, i.e. $E(n)_{p,q}$ where $p$ and $q$ are relatively prime. Later this method has proven useful in the construction of many other interesting 4-manifolds. This rational blow-down surgery was generalized by Jongil Park \cite{P1}. Below we discuss his generalized rational blow-down. Let $p\geq q\geq 1$ and $p,q$ are relatively prime. Let  $C_{p,q}$ be the smooth $4$-manifold obtained by plumbing disk bundles over the $2$-sphere according to the following linear diagram

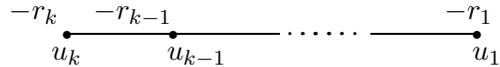
\begin{picture}(100,60)(-90,-25)
 \put(-12,3){\makebox(200,20)[bl]{$-r_{k}$ \hspace{6pt}
                                  $-r_{k-1}$ \hspace{96pt} $-r_{1}$}}
 \put(4,-25){\makebox(200,20)[tl]{$u_{k}$ \hspace{25pt}
                                  $u_{k-1}$ \hspace{86pt} $u_{1}$}}
  \multiput(10,0)(40,0){2}{\line(1,0){40}}
  \multiput(10,0)(40,0){2}{\circle*{3}}
  \multiput(100,0)(5,0){4}{\makebox(0,0){$\cdots$}}
  \put(125,0){\line(1,0){40}}
  \put(165,0){\circle*{3}}
\end{picture}

\noindent where $p^2/(pq-1) = [r_k, r_{k-1}, ..., r_1]$ is the unique continued linear fraction with all $r_{i} \geq 2$ and each vertex $u_{i}$ of the linear diagram represents a disk bundle over 2-sphere with Euler number $-r_{i}$. According to Casson and Harer \cite {CH}, the boundary of $C_{p,q}$ is the lens space $L(p^2, 1-pq)$ which also bounds a rational ball $B_{p,q}$ with $\pi_1(B_{p,q})=\Z_p$ and $\pi_1(\bd B_{p,q})\to \pi_1(B_{p,q})$ surjective. If $C_{p,q}$ is embedded in a $4$-manifold $X$ then the generalized rational blowdown manifold $X_{p,q}$  is obtained by replacing $C_{p,q}$ with $B_{p,q}$, i.e., $X_{p,q} = (X\- C_{p,q}) \cup B_{p,q}$. If $X$ and $X\- C_{p,q}$ are simply connected, then so is $X_{p,q}$. Notice that the case when $q=1$ is the  construction of Fintushel-Stern with $C_p = C_{p,1}$ given by

 
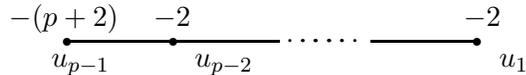
\begin{picture}(100,60)(-90,-25)
 \put(-12,3){\makebox(200,20)[bl]{$-(p+2)$ \hspace{6pt}
                                  $-2$ \hspace{96pt} $-2$}}
 \put(4,-25){\makebox(200,20)[tl]{$u_{p-1}$ \hspace{25pt}
                                  $u_{p-2}$ \hspace{86pt} $u_{1}$}}
  \multiput(10,0)(40,0){2}{\line(1,0){40}}
  \multiput(10,0)(40,0){2}{\circle*{3}}
  \multiput(100,0)(5,0){4}{\makebox(0,0){$\cdots$}}
  \put(125,0){\line(1,0){40}}
  \put(165,0){\circle*{3}}
\end{picture}

\smallskip 

\begin{lem} Let $X_{p,q}$ be the manifold obtained from $X$ by a rational blow-down of the configuration $C_{p,q}$ .  Then ${b_{2}}^{+}(X_{p,q}) = {b_{2}}^{+}(X)$ and ${c_{1}}^{2}(X_{p,q}) = {c_{1}}^{2}(X) + k$

\end{lem}

\begin{proof} 

Notice that the manifold $C_{p,q}$ is negative definite, so that ${b_{2}}^{+}(X_{p,q}) = {b_{2}}^{+}(X)$. Using the fact that ${c_1}^{2} = 3\s +2e$, we have ${c_{1}}^{2}(X_{p,q}) = 3\s(X_{p,q}) + 2e(X_{p,q}) = 3(\s(X)+k) + 2(e(X)-k) = {c_{1}}^{2}(X) + k$ 

\end{proof}

\begin{thm} \cite {P1}. Suppose $X$ is a smooth 4-manifold with $b_{2}^{+}(X) > 1$ which contains a configuration $C_{p,q}$. If $L$ is a characteristic line bundle on $X$ such that, $SW_{X}(L) \ne 0$, $(L|_{C_{p,q}})^{2} = - b_{2}(C_{p,q})$ and $c_{1}(L|_{L_(p^2,1-pq)}) = mp \in {\bf Z}_{p^{2}} \cong H^{2}(L(p^2, 1-pq); \bf Z)$ with $m \equiv (p-1) \mod 2$, then $L$ induces a SW basic class $\bar L$ of $X_{p,q}$ such that $SW_{X_{p,q}}(\bar L) = SW_{X}(L)$. \end{thm}

\medskip

\begin{coro} {\bf 3.2.} \cite {P3}. Suppose $X$ is a smooth 4-manifold with $b_{2}^{+}(X) > 1$ which contains a configuration $C_{p,q}$. If $L$ is a SW basic class of $X$ satisfying $L\cdot u_{i} = (r_{i} - 2)$ for any i with $1 \leq i \leq k$ (or  $L\cdot u_{i} = -(r_{i} - 2)$, then $L$ induces a SW basic class $\bar L$ of $X_{p,q}$ such that $SW_{X_{p,q}}(\bar L) = SW_{X}(L)$.  

\end{coro}

\section{Construction of exotic manifolds}

In this section, we construct a family of simply connected, nonsymplectic and nondiffemorphic manifolds starting from the elliptic surfaces $E(n)$ for $n\geq 3$ and applying the combination of knot surgery, blowup and rational blowdown. Using the Seiberg-Witten invariants, we distinguish their smooth structures. 
\begin{thm}

There are infinetly many simply connected, nonsymplectic and pairwise nondiffeomorphic manifolds that can be obtained from elliptic surfaces $E(n)$ $n\geq 3$ by the sequence of knot surgery in double nodes, ordinary blow ups and rational blowdowns with $(\chi_{h}, {c_{1}}^{2})$ given by following formulas.

\begin{displaymath}
\left\{\begin{array}{ll}   
\chi_{h} = 3, \ {c_{1}}^{2} \leq 16 \\
\chi_{h} = 4, \ {c_{1}}^{2} \leq 23 \\
\chi_{h} = 5, \ {c_{1}}^{2} \leq 30 \\
\chi_{h} = 6, \ {c_{1}}^{2} \leq 36 \\
\chi_{h} = n \geq 7, \ {c_{1}}^{2} \leq 25k - 2  & \textrm{if $n = 4k$}\\
\chi_{h} = n \geq 7, \ {c_{1}}^{2} \leq 25k + 5  & \textrm{if $n = 4k+1$}\\
\chi_{h} = n \geq 7, \ {c_{1}}^{2} \leq 25k + 11  & \textrm{if $n = 4k+2$}\\
\chi_{h} = n \geq 7, \ {c_{1}}^{2} \leq 25k + 18  & \textrm{if  $n = 4k+3$}\\
 
\end{array} \right.
\end{displaymath}

\end{thm}

\begin{proof}

Let first outline the general procedure how construct such family of simply connected, nonsymplectic 4-manifolds. According to lemma 1.1 $E(n)$ has elliptic fibration with $(2n+2)$ singular fibers, which are an $I_{8n}$ singularity, $(2n-1)$ number of $I_2$ singularity and two fishtail fibers. First perform a knot surgery in $s$ of the double nodes \cite{FS2} of the $(2n-1)$ $I_{2}$ singularity using non fibered twist knots $K_{1}, K_{2}, ..., K_{s}$, where $s = \lbrack \frac {7n+3}{4} \rbrack $. We denote the resulting manifold as $Y_{K_{1}, K_{2},.., K_{s}}$. The manifold  $Y_{K_{1},K_{2},.., K_{s}}$ is simply connected and has pseudo section $S$ \cite{FS2} which is sphere of self-intersection $-n$ with $s$ positive double points. The pseudo section $S$ intersects only one of the $-2$ spheres of $I_{8n}$ singularity. To construct our manifolds, we blow up $s$ double points on $S$ to get sphere $S'$ of self-intersection $-(n+4s)$ in $Y_{K_{1},K_{2},.., K_{s}}\#s\CPb$. Next blowup one of the fisthails if needed and smooth out its intersection with $S'$ to get sphere $S''$ of the self-intersection $-(8n+3)$. Then remove one of the spheres from the $I_{8n}$ configuration to get $(8n-1)$ linear chain of $-2$ spheres which together with $S''$ form the configuration $C_{8n+1}$. By blowing down of this configuration we get our manifolds. For simplicity let assume $K_{1} = K_{2} = K_{3} = .... = K_{r} = T(r)$ where $T(r)$ is non fibered $r$-twist knot. We will denote the manifold $Y_{K_{1},K_{2},..,K_{s}}$ in this set up as $Y(n)_{r, s}$ and the manifold obtained after rational blowdown as $Z(n)_{s, r, m}$, where $m$ denotes the number of the blowups applied. Notice that above construction is well defined once $n\geq 7$. The case when $n < 7$ can be treated similarly with small modification of the construction above and by rationally blowing down the configurations of smaller lenght.

Let first assume that $n \geq 7$. 

{\bf{Case 1}}. $n = 4k$. By Lemma 1.1, there is fibration with $I_{32k}$ singularity and $(2n-1) = 8k-1$ number of $I_{2}$ singularity and two fishtail fibers. Perform the knot surgery in $7k$ double nodes, blow up $s = 7k$ double points to get a sphere $S'$ with self-intersection $-32k$. Next blowup one of the fisthails, smooth out its intersection with $S'$ to get sphere $S''$ with self-intersection $-(32k+2)$. Use the linear chain of $(32k-2)$ of $-2$ spheres from $I_{32k}$ singularity together with $S''$ to get the configuration of $C_{32k}$ inside of $Y(n)_{r, s}\#(7k+1)\,\CPb$. By rationally blowing down this configuration we get manifolds with ${c_{1}}^{2} = 25k - 2$ $\chi_{h} = 4k$. 

{\bf{Case 2}}. $n = 4k+1$ By Lemma 1.1, there is fibration of $E(n)$ with one $I_{32k+8}$ singularity and $(2n-1) = 8k+1$ number of $I_{2}$ singularity and two fishtail fibers. After performing the knot surgery in $s = 7k +2$ double nodes, blow up the double points of pseudo-section, one fisthtail fiber and smooth out intersection of resulting sphere of self-intersection $-(32k+9)$ with fishtail fiber one gets the configuration of $C_{32k+9}$ inside of $Y(n)_{r, s}\#(7k+3)\,\CPb$. After rationally blowing down this configuration, we get manifolds with ${c_{1}}^{2} = 25k + 5$ and $\chi_{h} = 4k + 1$. 

{\bf{Case 3}}. $n = 4k+2$. By Lemma 1.1, there is fibration with $I_{32l+16}$ singularity and $(2n-1) = 8k+3$ number of $I_{2}$ singularity and two fishtail fibers. Perform knot surgery in $7k+4$ double nodes, blow up $s = 7k +4$ times to get a sphere $S'$ with self-intersection $-(32k+18)$. Next use the linear chain of $(32k+14)$ of $-2$ spheres from $I_{32k+16}$ singularity together with $S'$ to get configuration of $C_{32k+16}$ inside of $Y(n)_{r, s}\#(7k+4)\,\CPb$. Rational blowdown of this configuration gives the manifolds with ${c_{1}}^{2} = 25k + 11$ $\chi_{h} = 4k + 2$ 

{\bf{Case 4}}. $n = 4k+3$. Again by Lemma 1.1, $E(n)$ has fibration with $I_{32k+24}$ singularity and $(2n-1) = 8k+5$ number of $I_{2}$ singularity and two fishtails. We perform knot surgery on $s = 7k + 6$ double nodes and blow up double points on pseudo-section to get configuration of $C_{32k+25} = C_{8n+1}$ configuration inside of $Y(n)_{r, s}\#(7k+6)\,\CPb$. By rational blowdown we get manifold with ${c_{1}}^{2} = 25k + 18$ and $\chi_{h} = 4k + 3$ 

\smallskip

Notice that these manifolds are simply-connected due to the fact that there is a fishtail fiber that was not used in the construction. By computing the Seiberg-Witten invariants, we will distinguish the different smooth structures constructed. 

\smallskip

To compute Seiberg-Witten invariants of the manifold $Y(n)_{r, s}$, one applies the result of Fintushel-Stern on how Seiberg-Witten invariants change under the knot surgery \cite {FS2} . Using the facts that the Seiberg-Witten function of $E(n)$ is $(e^{T} - e^{-T})^{n-2}$ and the Alexander polynomial of the twist knot $T(r)$ is $\DD_{T(r)} = (2r-1) + rt - rt^{-1}$, one computes the Seiberg-Witten invariants of $SW_{Y_{n,r}} = (e^{T} - e^{-T})^{n-2} ((2r-1) + re^{2T} - re^{-2T})^{s}$. By blow up formula for the Seiberg-Witten function \cite{FS2}, we have $SW_{Y(n)_{r, s}\#m\,\CPb} = SW_{Y_{n,r}} \cdot \prod_{j=1}^{m}(e^{E_{i}} + e^{-E_{i}})$, where $E_{i}$ is an exceptional class coming from the $i^{th}$ blow up. Furthermore, using Corollary 3.2 and Theorem 1.2 \cite {P1},  we compute the Seiberg-Witten invariants of $Z(n)_{s,r,m}$: Up to sign it has only one basic class which corresponds to top classes ${\pm ((n+ 2s - 2)T + E_{1} + ... + E_{m})}$ and the value of the Seiberg-Witten function on these classes are ${ \pm r^{s}}$. Now by applying the theorem of Taubes from  Section 2, we detect that these manifolds are not symplectic. Since the Seiberg-Witten invariants are diffeomorphism invariants, the manifolds constructed yield an infinite family of simply connected, nonsymplectic and pairwise nondiffeomorphic manifolds once we consider different non fibered twist knots in the construction.  

\medskip

{\bf{Case $n<7$}}. We carry out the computation for the cases $n = 3, 4$. The remaning cases are similar. 

{\bf{i. $n=3.$}} According to Lemma 1.1, $E(3)$ have one $I_{24}$ singularity, $5$ $I_{2}$ singularity and two fishtail fibers. Perform knot surgery in all 5 double nodes, blowup $5$ double points of the pseudo section $S$, one of the fishtails and smooth out its intersection with the pseudo section to get sphere $S'$ of self interesction $-25$. By using $21$ of the $-2$ spheres of the $I_{24}$ singularity and $S'$ one gets the configuration of $C_{23}$ inside of $E(3)\#6\,\CPb$. After rational blowdown of $C_{23}$, one gets simply-connected manifolds with ${c_{1}}^{2} = 16$, $\chi_{h} = 3$. 

{\bf{ii. $n=4.$}} Again, by Lemma 1.1, $E(4)$ admits a fibration with one $I_{32}$ singularity, $7$ $I_{2}$ singularity and two fishtail fibers. Perform knot surgery in all $7$ double nodes, blowup $7$ double points of the pseudo section $S$, one of the fishtails and smooth out its intersection with the pseudo section to get sphere $S'$ of self interesction $-34$. By using $30$ of the $-2$ spheres of the $I_{32}$ singularity and $S'$ one gets the configuration of $C_{32}$ inside of $E(4)\#8\,\CPb$. After rational blowdown of $C_{32}$, we get simply-connected manifolds with ${c_{1}}^{2} = 23$, $\chi_{h} = 4$.

 The computation of Seiberg-Witten invariants is similar. 

Now varying the number of double nodes $s$ from $0$ to $\lbrack {7n+3}/4 \rbrack$ and applying the blow ups, we can fill all the lattice points given as in the statement of the theorem. Notice that the number of basic classes will increase and ${c_{1}}^{2}$ reduce by blow up operation. 

\end{proof}

\end{document}